\documentclass[a4paper, 11pt, oneside, final]{scrartcl}

\usepackage{PREAMBLE}
\usepackage{COMMAND}
\usepackage{HYPHENATION}

\begin{document}
\thispagestyle{empty}
\begin{center}
    \rule{\linewidth}{1pt}\\[0.4cm]
    {\sffamily \bfseries \large One-dimensional phase retrieval with\\ additional
        interference measurements}\\[10pt] 
    {\sffamily\footnotesize Robert Beinert}\\[3pt]
    {\sffamily\footnotesize Institut für Numerische und Angewandte Mathematik}\\[-3pt]
    {\sffamily\footnotesize Georg-August-Universität Göttingen}\\
    \rule{\linewidth}{1pt}
\end{center}

\vspace*{10pt}

\begin{center}
    \sffamily\bfseries Abstract
\end{center}

{\small
    \noindent
    The one-dimensional phase retrieval problem consists in the recovery of a
    complex-valued signal from its Fourier intensity.  Due to the well-known ambiguousness
    of this problem, the determination of the original signal within the extensive
    solution set is challenging and can only be done under suitable a~priori assumption or
    additional information about the unknown signal.  Depending on the application, one
    has sometimes access to further interference measurements between the unknown signal
    and a reference signal.  Beginning with the reconstruction in the discrete-time
    setting, we show that each signal can be uniquely recovered from its Fourier intensity
    and two further interference measurements between the unknown signal and a modulation
    of the signal itself.  Afterwards, we consider the continuous-time problem, where we
    obtain an equivalent result.  Moreover, the unique recovery of a continuous-time
    signal can also be ensured by using interference measurements with a known or an
    unknown reference which is unrelated to the unknown signal.

    \medskip

    \noindent
    {\sffamily\bfseries Key words:} one-dimensional phase retrieval for complex
    signals, compact support, interference measurements

    \smallskip
    
    \noindent
    {\sffamily\bfseries AMS Subject classifications:} 42A05, 94A08, 94A12, 94A20
}

\section{Introduction}

Phase retrieval problems occur in a width range of applications in physics and engineering
such as crystallography \cite{Mil90,Hau91,KH91}, astronomy \cite{BS79,DF87}, and laser
optics \cite{SST04,SSD+06}.  All of these applications have in common that one needs to
recover an unknown signal from the intensity of its Fourier transform.  Because of the
well-known ambiguousness of this problem, the recovery of an analytic or a numerical
solution is generally challenging.  To determine a meaningful solution within the solution
set, one therefore needs further a priori information or additional data.

In this paper, we consider the one-dimensional phase retrieval problem for discrete-time
and continuous-time signals.  In both cases, we restrict ourselves to the recovery of an
unknown signal with finite or compact support.  The ambiguities of these problems can be
explicitly specified by an appropriate factorization of the corresponding autocorrelation
signal, see \cite{BS79,OS89,BP15} for the discrete-time and \cite{Aku56,Aku57,Wal63,Hof64}
for the continuous-time setting.

Depending on the application, one can sometimes superpose the unknown signal $x$ with an
appropriate reference signal $h$.  In the discrete-time setting, the additional Fourier
intensity of the interference $x+h$ can be used to reduce the solution set to merely two
different signals or even to ensure uniqueness for particularly known references,
\cite{KH90a,KH90,BP15}.  Under some further assumptions, it is also possible to use
unknown references to recover the wanted signal, \cite{KH93,RDN13,RSA+11,BP15}.  Besides
considering reference signals being not related to the unknown signal, one may also employ
a modulation of the unknown signal itself as reference in order to guarantee uniqueness,
\cite{CESV13}.

A slightly different approach to ensure uniqueness by additional data in the frequency
domain is presented in \cite{NQL83a,NQL83}.  Here the Fourier transform is replaced by the
so-called short-time Fourier transform, where the unknown signal is overlapped with a
small analysis window at different positions.

If we normalize the support of the unknown signal $x$ to $\{0,\dots,M-1\}$, we can
interpret $x$ as an $M$-dimensional vector.  Further, the Fourier intensities
$\absn{\fourier x(\omega_k)}$ at different points $\omega_k \in [-\pi,\pi)$ can now be
written as the intensity measurement $\absn{\iProd{x}{v_{k}}}$ with
$v_{k} \coloneqq (\e^{\I \omega_k m})_{m=0}^{M-1}$.  Generalizing this idea, here the
question arises how the vectors $v_k$ have to be chosen, and how many vectors $v_k$ are
required to ensure the recovery of $x$ from the intensities $\absn{\iProd{x}{v_k}}$, see
for instance \cite{BCE06, BBCE09,BCM14,BH15} and references therein.  To reconstruct the
signal $x$ explicitly, one can again exploit suitable interference measurements of the
form $\absn{\iProd{x}{v_k+v_\ell}}$, see \cite{ABFM14}.

Besides employing interference measurements, the ambiguities can also be avoided by
further information about the unknown signal in the time domain.  In \cite{BP16}, it has
been shown that almost every complex-valued discrete-time signal $x$ is uniquely
determined by its Fourier intensity $\absn{\fourier x}$ and one absolute value
$\absn{x[n]}$ for a fixed $n$ within the support of $x$.  A similar statement holds if two
phases $\arg x[n]$ and $\arg x[m]$ for appropriate $n$ and $m$ are known beforehand.

The continuous-time phase retrieval problem, however, has a completely different behavior
than the discrete-time equivalent.  Nonetheless, also here the ambiguousness is a
challenging problem.  If the additional Fourier intensity of an appropriate modulation of
the unknown signal is available, one can here solve the corresponding phase retrieval
problem uniquely by comparing the zeros of the analytic continuation of the given
intensities, \cite{WFB81}.  Using a combination of oversampling and modulations, Pohl
\etal\ show that the unknown signal can be recovered up to a global phase, \cite{PYB14a}.
Similarly to the discrete-time setting, the unknown signal can also be uniquely recovered
by exploiting interference measurements.  For this purpose, Burge \etal\ have constructed
suitable reference signals depending on the given Fourier intensity, \cite{BFGR76}.

The paper is organized as follows.  In section~\ref{sec:discrete-time-phase}, we introduce the
discrete-time phase retrieval problem.  Here we investigate interference measurements of
the unknown signal with a modulation of the signal itself.  Using Prony's method, we adapt
and extend the uniqueness results in \cite{ABFM14} and \cite{CESV13} to our specific
problem formulation and show that each discrete-time signal with finite support can always
be recovered uniquely.  Moreover, we can drop the assumption that the support of the
unknown signal must be known beforehand as supposed in \cite{ABFM14} and \cite{CESV13},
since we can recover the support from the given Fourier intensities.

In section~\ref{sec:cont-time-phase}, we transfer our findings from the discrete-time to the
continuous-time setting.  Moreover, we consider interference measurements where the
reference signal is known beforehand.  Similar to the discrete-time setting, these
measurements reduce the solution set to two different signals.  Giving a novel proof of
the corresponding results in \cite{KST95}, we can here explicitly represent the second
possible solution.  Moreover, we show that the uniqueness results for an unknown reference
signal \cite{KH93,RDN13,BP15} can also be transferred.  For this purpose, we will employ
the characterization of the arising ambiguities in \cite{Hof64}.

\section{Discrete-time phase retrieval}
\label{sec:discrete-time-phase}

\subsection{Trivial and non-trivial ambiguities}
\label{sec:trivial-non-trivial}

We begin with the \emph{one-dimensional discrete-time phase retrieval problem}, where we
wish to recover a complex-valued signal $x \coloneqq (x[n])_{n\in\Z}$ from its Fourier
intensity $\absn{\fourier x}$, where $\fourier x$ denotes the \emph{discrete-time Fourier
    transform} given by
\begin{equation*}
    \fourier x \mleft( \omega \mright)
    \coloneqq \Fourier \mleft[ x \mright] \mleft( \omega \mright)
    \coloneqq \sum_{n\in\Z} x\mleft[n\mright] \, \e^{-\I \omega n}
    \qquad (\omega \in \R).
\end{equation*}
In the following, we assume that the unknown signal $x$ has a finite support of length
$N$.  In other words, we always find an $n_0 \in \Z$ such that $x[n] = 0$ for $n < n_0$
and $n \ge n_0 +N$.  As a consequence of this a priori condition, the corresponding
\emph{autocorrelation signal}
\begin{equation*}
    a \mleft[ n \mright] \coloneqq
    \sum_{k\in\Z} \overline{x \mleft[ k \mright]} \, x \mleft[ k+n \mright]
    \qquad
    (n \in \Z)
\end{equation*}
possesses the finite support $\{-N+1, \dots, N-1 \}$.  Moreover, the
\emph{autocorrelation function}
\begin{equation*}
    \fourier a \mleft( \omega \mright)
    = \sum_{n\in \Z} a \mleft[ n \mright] \, \e^{-\I\omega n}
    = \sum_{n\in \Z} \sum_{k\in\Z} x \mleft[n\mright] \, \overline{x \mleft[k\mright]} \,
    \e^{-\I \omega (n-k)}
    =\absn{\fourier x \mleft( \omega \mright)}^2
\end{equation*}
is here always a non-negative trigonometric polynomial of degree $N-1$, and the support
length $N$ of the unknown signal $x$ is hence explicitly encoded in the given Fourier
intensity $\absn{\fourier x}$.

Under the assumption that the support length $N$ is known beforehand, the considered phase
retrieval problem can be discretized in the frequency domain.  Our
problem is equivalent to the task to recover the unknown signal $x$ with support length
$N$ from the $2N-1$ values
\begin{equation*}
    \abs{\fourier x \mleft( \tfrac{2\pi k}{N} \mright)}
    \qquad
    (k=-N,\dots,N-1)
    \addmathskip
\end{equation*}
since the trigonometric polynomial $\fourier a = \absn{\fourier x}^2$ of degree $N-1$ is completely
determined by $2N-1$ samples at different points in $[-\pi, \pi)$.

It is well known that the (discrete-time) phase retrieval problem is not uniquely solvable
in general.  The simplest occurring ambiguities are caused by rotation (multiplication
with a unimodular constant), time shift, and conjugation and reflection, \cite{BP15}.

\begin{Proposition}\label{prop:trivial-amb}
    Let $x$ be a complex-valued signal with finite support. Then
    \begin{enumerate}[\upshape(i)]
    \item\label{item:1} the rotated signal $( \e^{\I \alpha} \, x[ n] )_{n\in\Z}$ for
        $\alpha \in \R$
    \item\label{item:2} the time shifted signal $( x[n-n_0] )_{n\in\Z}$ for $n_0 \in \Z$
    \item\label{item:3} the conjugated and reflected signal $(\overline{ x[-n]})_{n\in\Z}$
    \end{enumerate}
    have the same Fourier intensity $\absn{\fourier x}$.%
\end{Proposition}

Although the ambiguities in Proposition~\ref{prop:trivial-amb} always occur, they are of minor
interest since they are closely related to the original signal.  Usually, besides the
rotations, shifts, and conjugation and reflection, the discrete-time phase retrieval
problem has an extensive set of further solutions, which can completely differ from the
original signal.  Similarly to \cite{Wan13,BP15}, we distinguish between trivial and
non-trivial ambiguities.

\begin{Definition}
    \label{def:tri-nontri-amb}%
    A \emph{trivial ambiguity} of the discrete-time phase retrieval problem is caused by a
    rotation, time shift, conjugation and reflection, or by a combinations of these. All
    other occurring ambiguities are called \emph{non-trivial}.%
\end{Definition}

Unlike the trivial ambiguities, where we have infinitely many possibilities to choose the
rotation parameter $\alpha \in \R$ and the shift parameter $n_0 \in \Z$, the phase
retrieval problem to recover a signal with finite support can only possesses finitely many
non-trivial solutions, see \cite[Corollary~2.6]{BP15}.  More precisely, depending on the
support length $N$ of the unknown signal, the solution set can consists of at most
$2^{N-2}$ non-trivially different signals.

\subsection{Interference measurements in the discrete-time setting}
\label{sec:interf-modul-sig}

Considering the large number of non-trivial solutions, how can we determine the original
signal within the set of ambiguities?  One possibility to reduce the solution set
significantly is the exploitation of additional interference measurements of the form
$\absn{\Fourier [ x + h ]}$, where $h$ is a suitable reference signal with finite support.
If $h$ is a known reference signal unrelated to $x$, then the corresponding phase
retrieval problem has at most two non-trivially different solutions, see \cite{KH90,BP15}.
If $h$ is also an unknown signal, under some additional assumption, both signals $x$ and
$h$ can be uniquely recovered from $\absn{\fourier x}$, $\absn{\fourier h}$, and
$\absn{\fourier x + \fourier h}$, see \cite{KH93,RDN13,BP15}, where in \cite{RDN13} the
additional measurements $\absn{\fourier x + \I \fourier h}$ are employed.

In this section, we examine interference measurements of a slightly different kind.  More
precisely, we replace the reference signal $h$ with a modulated version of the unknown
signal $x$ itself.  This idea goes back to \cite{CESV13} and \cite{ABFM14}, where the
phase retrieval of a finite-dimensional vector from the intensities of the discrete
Fourier transform or the intensities of a suitably constructed frame is considered.

Transferring this approach to the discrete-time phase retrieval problem, we try to recover
the complex-valued signal $x$ with finite support from its Fourier intensity
$\absn{\fourier x}$ and a set of interference measurements of the form
\begin{equation}
    \label{eq:int-meas}
    \abs{\Fourier \mleft[ x + \e^{\I \alpha} \, \e^{\I \mu \cdot} \, x\mright]} \!,
\end{equation}
where the modulations and rotations are described by $\mu\in \R$ and $\alpha \in \R$.
Since the Fourier transform of the modulated signal is given by
\begin{equation*}
    \Fourier \mleft[ \e^{\I \alpha} \, \e^{\I \mu \cdot} \, x  \mright] \mleft( \omega
    \mright) 
    = \e^{\I \alpha} \sum_{n \in \Z} x \mleft[ n \mright]  \, \e^{- \I \left( \omega - \mu \right) n} 
    = \e^{\I \alpha} \, \fourier x \mleft( \omega - \mu \mright),
\end{equation*}
we can also interpret the considered interference measurements as interferences with
certain shifts of the Fourier transform $\fourier x$ in the frequency domain.

\subsubsection{Phase reconstruction by using a polarization identity}

Our first approach to recover $x$ is to apply a suitable polarization identity, which
allows us to determine the unknown phase of $\fourier x$.  For this purpose, we generalize
the Mercedes-Benz polarization identity introduced by Alexeev \etal\ in \cite{ABFM14}.  In
the following, the primitive $K$th root of unity is denoted by
$\zeta_K \coloneqq \e^{\nicefrac{2\pi\I}{K}}$.

\begin{Lemma}
    \label{lem:polar-id}%
    Let $z_1$ and $z_2$ be two complex numbers.  Then the polarization identity
    \begin{equation}
        \label{eq:polar-id}
        \overline z_1 z_2
        = \frac{1}{K} \sum_{k=0}^{K-1} \zeta_K^k \abs{ z_1 + \zeta_K^{-k} \, z_2}^2
    \end{equation}
    holds for every integer $K>2$.%
\end{Lemma}

\begin{Proof}
    The assertion can be proven by generalizing the ideas in \cite[\p~38]{ABFM14}.  We expand the
    right-hand side of \eqref{eq:polar-id} and obtain
    \begin{align*}
      \frac{1}{K} \sum_{k=0}^{K-1} \zeta_K^k \abs{ z_1 + \zeta_K^{-k} \, z_2}^2
      &= \frac{1}{K} \sum_{k=0}^{K-1} \zeta_K^k \left( \abs{z_1}^2 + 2
        \Re \mleft[ \zeta_K^{-k} \, \overline z_1 z_2 \mright]
        + \abs{z_2}^2 \right)
      \\[\fskip]
      &= \frac{2}{K} \sum_{k=0}^{K-1} \zeta_K^k \, \Re \mleft[ \zeta_K^{-k} \, \overline z_1
        z_2 \mright]
    \end{align*}
    since the sum over all roots $\zeta_K^k$ is zero.  Writing the real part of the
    product $\zeta_K^{-k} \, \overline z_1 z_2$ as
    \begin{equation*}
        \Re \bigl[ \zeta_K^{-k} \, \overline z_1  z_2 \bigr]
        = \Re \bigl[ \zeta_K^{-k} \bigr] \, \Re \bigl[ \overline z_1  z_2 \bigr]
        - \Im \bigl[ \zeta_K^{-k} \bigr] \, \Im \bigl[ \overline z_1  z_2 \bigr],
    \end{equation*}
    and using the identities $ \Re \zeta_K^{-1} = \Re \zeta_K^k$ and
    $ \Im \zeta_K^{-k} = - \Im \zeta_K^k$, we further have
    \begin{equation*}
        \frac{1}{K} \sum_{k=0}^{K-1} \zeta_K^k \abs{ z_1 + \zeta_K^{-k} \, z_2}^2
        = \frac{2}{K} \sum_{k=0}^{K-1} \zeta_K^k 
        \left( \Re \bigl[ \zeta_K^{k} \bigr] \, \Re \bigl[ \overline z_1  z_2 \bigr]
            + \Im \bigl[ \zeta_K^{k} \bigr] \, \Im \bigl[ \overline z_1  z_2 \bigr]
        \right) \!.
    \end{equation*}

    We consider the real and imaginary part of this equation separately.  Here the
    substitution $\zeta_K^k \coloneqq \Re \zeta_K^k + \I \Im \zeta_K^k$ yields
    \begin{align*}
      \Re \mleft[\frac{1}{K} \sum_{k=0}^{K-1} \zeta_K^k \abs{ z_1 + \zeta_K^{-k} \, z_2}^2 \mright]
      &= \frac{2}{K} \left( \Re \bigl[ \overline z_1  z_2 \bigr] \sum_{k=0}^{K-1}
        \bigl[ \Re \zeta_K^{k} \bigr]^2
        + \Im \bigl[ \overline z_1 z_2 \bigr] \sum_{k=0}^{K-1}
        \Re \zeta_K^{k} \, \Im \zeta_K^{k} \right)
        \intertext{and}
        \Im \mleft[\frac{1}{K} \sum_{k=0}^{K-1} \zeta_K^k \abs{ z_1 + \zeta_K^{-k} \, z_2}^2 \mright]
      &= \frac{2}{K} \left( \Re \bigl[ \overline z_1 z_2 \bigr] \sum_{k=0}^{K-1}
        \Im \zeta_K^k \, \Re \zeta_K^{k} 
        + \Im \bigl[ \overline z_1  z_2 \bigr] \sum_{k=0}^{K-1}
        \bigl[  \Im \zeta_K^{k} \bigr]^2  \right) \!.
    \end{align*}
    Rewriting the real and imaginary parts as
    $\Re \zeta_K^k = \nicefrac{1}{2} (\zeta_K^{-k} + \zeta_K^{k})$ and
    $\Im \zeta_K^k = \nicefrac{\I}{2} (\zeta_K^{-k} - \zeta_K^{k})$, we can compute the
    occurring sums by
    \begin{equation*}
        \sum_{k=0}^{K-1} \bigl[\Re \zeta_K^k \bigr]^2
        = \frac{1}{4} \sum_{k=0}^{K-1} \left( \zeta_K^{-2k} + 2 + 
            \zeta_K^{2k} \right)
        = \frac{K}{2},
        \addmathskip
    \end{equation*}
    \begin{equation*}
        \sum_{k=0}^{K-1} \bigl[\Im \zeta_K^k \bigr]^2
        = -\frac{1}{4} \sum_{k=0}^{K-1}\left( \zeta_K^{-2k} - 2 + 
            \zeta_K^{2k} \right)
        = \frac{K}{2},
    \end{equation*}
    and 
    \begin{equation*}
        \sum_{k=0}^{K-1} \Re \zeta_K^k \Im \zeta_K^k
        = \frac{\I}{4} \sum_{k=0}^{K-1} \left( \zeta_K^{-2k} - \zeta_K^{2k} \right)
        = 0,
        \addmathskip
    \end{equation*}
    which completes the proof. \qed%
\end{Proof}

\begin{Remark}
    \label{rem:polar-id}
    Obviously, Lemma~\ref{lem:polar-id} cannot be valid if $K<3$ because the right-hand side
    would be always a real number.  For the special case $K=3$, Lemma~\ref{lem:polar-id}
    coincides with the original Mercedes-Benz polarization identity introduced in
    \cite[Lemma~2.1]{ABFM14}.  \qed%
\end{Remark}

The polarization identity in Lemma~\ref{lem:polar-id} can now be used
to reveal the relation between the values $\fourier x(\omega)$ and
$\fourier x (\omega - \mu)$ hidden in the interference measurements
\begin{equation*}
    \abs{\Fourier \mleft[ x + \zeta_K^{-k} \, \e^{\I \mu \cdot} \,  x \mright]}
    \qquad \left( k = 0, \dots, K-1 \right)
\end{equation*}  
for some $K \ge 3$ and $\mu \in \R$.  The following theorem shows that the knowledge of
this relationship is sufficient to recover each discrete-time signal up to trivial
rotations.

\begin{Theorem}
    \label{the:interf-polar-id}%
    Let $x$ be a discrete-time signal with finite support of length $N$.  If
    $\mu \ne 2 \pi \,\nicefrac{p}{q}$ for all $p\in\Z$ and $q \in \{1,\dots,N-1\}$, then
    the signal $x$ can be uniquely recovered up to rotation from its Fourier
    intensity $\absn{\fourier x}$ and the interference measurements
    \begin{equation*}
        \abs{\Fourier \mleft[ x + \zeta_K^{-k} \, \e^{\I \mu \cdot} \,  x \mright]}
        \qquad \left( k = 0, \dots, K-1 \right)
        \submathskip
    \end{equation*}   
    for every $K \ge 3$.
\end{Theorem}

\begin{Proof}
    Firstly, we apply the polarization identity in Lemma~\ref{lem:polar-id} to the given
    interference measurements
    $\absn{\fourier x ( \cdot ) + \zeta_K^{-k} \, \fourier x ( \cdot - \mu )}$
    and obtain
    \begin{equation*}
        \frac{1}{K} \sum_{k=0}^{K-1} \zeta_K^k \abs{\fourier x \mleft( \omega \mright) +
            \zeta_K^{-k} \, \fourier x \mleft( \omega - \mu\mright)}^2
        = \abs{\fourier x \mleft( \omega \mright)}
        \abs{ \fourier x \mleft( \omega - \mu \mright)}
        \e^{\I \left( \phi \mleft( \omega -\mu \mright) - \phi \mleft( \omega \mright) \right)},
    \end{equation*}
    where $\phi$ denotes the phase of $\fourier x = \absn{\fourier x} \, \e^{\I\phi}$.  If
    $\fourier x (\omega)$ and $\fourier x(\omega - \mu)$ are non-zero, we can hence
    determine the phase difference $\phi(\omega - \mu) - \phi( \omega)$.

    Since the unknown signal can only be recovered up to rotations, we can arbitrarily
    choose the phase $\phi(\omega_0)$ of $\fourier x (\omega_0)$ for one $\omega_0$
    without loss of generality.  Starting from this point, we can now use the phase
    differences to compute a series of relative phases
    \begin{equation*}
        \phi \mleft( \omega_0 + \mu k \mright)
        \qquad (k=0,\dots,2N-1)
    \end{equation*}
    with respect to $\phi(\omega_0)$.  To ensure that $\phi(\omega_0 + \mu k) \ne 0$ for
    $k=0,\dots, 2N-1$, we notice that $\fourier x$ can be written as
    \begin{equation}
        \label{eq:interf-polar-id:xhat}
        \fourier x \mleft( \omega \mright)
        = \e^{-\I \omega n_0} \sum_{n=0}^{N-1} c_n \, \e^{-\I \omega n},
        \addmathskip
    \end{equation}
    with $c_n \coloneqq x [n+n_0]$ because of the finite support.  Since the occurring
    algebraic polynomial in $\e^{-\I\omega}$ of degree $N-1$ can only have finitely many
    zeros on the unit circle, we can always find a suitable $\omega_0$ that enables us to
    recover $2N$ relative phases and hence $2N$ values of $\fourier x$ itself.  Writing
    the recovered points of $\fourier x$ in \eqref{eq:interf-polar-id:xhat} as
    \begin{equation}
        \label{eq:interf-polar-id:exp-sum}
        \fourier x \mleft( \omega_0 + \mu k\mright)
        = \sum_{n=0}^{N-1} \left[ c_n \, \e^{-\I \omega_0 (n+n_0)} \right]
        \e^{-\I k \mu (n+n_0)}
        \qquad(k=0, \dots, 2N-1),
    \end{equation}
    we can interpret the determined points as function values of the exponential sum with
    complex coefficients $c_n \, \e^{-\I \omega_0 (n + n_0)}$ and real frequencies
    $\mu(n+n_0)$ at the equally spaced points $k$ from $0$ to $2N-1$.

    By assumption, we can suppose that $\mu \ne 2 \pi \, \nicefrac{\ell}{(n_1-n_2)}$ for
    all $\ell \in \Z$ and $n_1 \ne n_2$ with $n_1$, $n_2 \in \{0, \dots, N-1 \}$. Therefore,
    \begin{equation*}
        \mu \left( n_1 + n_0 \right) \ne \mu \left( n_2 + n_0 \right) + 2\pi \ell
    \end{equation*}
    and thus the values $\e^{-\I \mu (n + n_0)}$ in \eqref{eq:interf-polar-id:exp-sum}
    differ pairwise for $n=0,\dots,N-1$.  Consequently, the exponential sum in
    \eqref{eq:interf-polar-id:exp-sum} is completely determined by the given function
    values, and we can apply Prony's method to determine the unknown coefficients and
    frequencies, see \cite{Pro95} or \cite[Sect.~9.4]{Hil87}.  For that reason, we can
    always determine $\fourier x$ and hence $x$ up to rotation. \qed%
\end{Proof}

\begin{Remark}
    \label{rem:interf-polar-id:1}
    The main reason for the application of Prony's method in the proof of
    Theorem~\ref{the:interf-polar-id} is the lack of information about the integer $n_0$ in the
    frequency representation \eqref{eq:interf-polar-id:xhat}.  Considering the influence
    of the modulation $\e^{-\I \omega n_0}$ in the time domain, we only know the support
    length $N$ of the unknown signal $x$ but not the exact position of the support itself.
    If we additionally assume that $x$ has the support $\{0, \dots, N-1 \}$, we can
    recover the trigonometric polynomial $\fourier x$ directly from the constructed
    function values by solving a linear equation system because all occurring frequencies
    are known beforehand.  Unfortunately, we cannot neglect the restrictions on the
    parameter $\mu$ since these are needed to ensure the invertibility of the arising
    Vandermonde matrix.  \qed%
\end{Remark}

\begin{Remark}
    \label{rem:interf-polar-id:2}
    Considering the assumptions of Theorem~\ref{the:interf-polar-id}, we have to choose $\mu$ in
    a way that $\mu$ is not a rational multiple of $2\pi$ where the denominator is an
    integer between $1$ and $N-1$.  Choosing $\mu$ as an irrational multiple of $2 \pi$,
    we can recover every signal from the given interference measurements independently of
    the actual support length $N$.  \qed%
\end{Remark}

\begin{Remark}
    \label{rem:interf-polar-id:3}
    If we assume that the support of the unknown discrete-time signal $x$ with support
    length $N$ is contained in the set $\{0,\dots, M-1\}$ with $M \ge N$ such that
    $(x[n])_{n\in\Z}$ can be identified with an $M$-dimensional vector, then we can
    interpret the Fourier intensity $\absn{\fourier x (\omega_n)}$ for a certain point
    $\omega_n$ in the frequency domain as intensity measurement $\absn{\iProd{x}{v_n}}$
    with the frame vector $v_n \coloneqq (\e^{\I\omega_n k} )_{k=0}^{M-1}$.  Choosing at
    least $M$ pairwise different points $\omega_n$ in $[-\pi, \pi)$ beforehand, we can
    consequently transfer the whole theory developed by Alexeev \etal{} \cite{ABFM14} to recover the
    unknown vector $x$ from the given intensity measurements
    $\absn{\fourier x (\omega_n)}$ and the given interference measurements
    \begin{equation*}
        \abs{\fourier x\mleft(\omega_n\mright)
            + \zeta_3^{-k} \, \fourier x \mleft(\omega_m\mright)}
        \qquad
        \left(k = 0, 1, 2 \right)
        \addmathskip
    \end{equation*}
    for a larger number of randomly chosen index pairs $(n,m)$.  In \cite{ABFM14} it has
    been shown that the unknown vector can be recovered with high probability from
    approximately $240 M$ measurements for an arbitrary frame, see \cite[\p~41]{ABFM14}.

    In contrast to the findings of Alexeev \etal{} for arbitrary frames, we exploit that
    the Fourier transform of a complex-valued discrete-time signal with finite support of
    length $N$ is a trigonometric polynomial.  For the special case $K=3$ considered by
    Alexeev \etal{}, this enables us to recover the unknown signal $x$ always from merely
    $8N-4$ measurements or, more precisely, from $2N-1$ measurements for each of the four
    Fourier intensities
    \begin{equation*}
        \abs{\fourier x\mleft(\cdot\mright) }
        \quad\text{and}\quad
        \abs{\fourier x\mleft(\cdot\mright) + \zeta_3^{-k} \, \fourier x \mleft( \cdot -
            \mu \mright)}
        \qquad (k=0,1,2),
        \addmathskip
    \end{equation*}
    \cf\ section~\ref{sec:trivial-non-trivial}.  Moreover, the procedure in the proof of
    Theorem~\ref{the:interf-polar-id} allows us to determine the position of the current support
    from the given Fourier intensities.  In other words, we can recover the unknown
    signal $x$ without the assumption that the support of $x$ is contained in some
    specific set $\{0,\dots, M-1\}$.  \qed%
\end{Remark}

\subsubsection{Reducing the number of required interference measurements}

Looking back at Theorem~\ref{the:interf-polar-id}, we observe that the actual number of
interference measurements depends on the chosen root of unity $\zeta_K$ or, more
precisely, on the chosen integer $K$.  Consequently, it seems that the given interference
measurements are highly redundant.  This impression is confirmed by a result in
\cite{CESV13}, where Cand\`es \etal\ employ only two of the interference measurements like
in Theorem~\ref{the:interf-polar-id} to recover a finite-dimensional vector.

In this section, we adapt the approach of Cand\`es \etal\ to the discrete-time phase
retrieval problem and simultaneously generalize the result to show that each
complex-valued discrete-time signal can be recovered from its Fourier intensity and two
further interference measurements of the form \eqref{eq:int-meas}, where the two rotations
can be chosen almost arbitrarily.  In particular, these rotations do not have to arise
from the $K$th roots of unity for a certain integer $K$.

\begin{Theorem}
    \label{the:two-interf-arb-rot}%
    Let $x$ be a discrete-time signal with finite support of length $N$. If
    $\mu \ne 2 \pi \, \nicefrac{p}{q}$ for all $p\in \Z$ and $q \in \{1,\dots,N-1\}$,
    then the signal $x$ can be uniquely recovered up to a rotation from its Fourier
    intensity $\absn{\fourier x}$ and the two interference measurements
    \begin{equation*}
        \abs{\Fourier \mleft[ x + \e^{\I \alpha_1} \, \e^{\I \mu \cdot} \, x\mright]}
        \quad\text{and}\quad
        \abs{\Fourier \mleft[ x + \e^{\I \alpha_2} \, \e^{\I \mu \cdot} \, x\mright]} \!,
    \end{equation*}
    where $\alpha_1$ and $\alpha_2$ are two real numbers satisfying
    $\alpha_1 - \alpha_2 \ne \pi k$ for all integer $k$.%
\end{Theorem}

\begin{Proof}
    We follow the lines of the proof of Theorem~\ref{the:interf-polar-id}.  Again, the crucial
    point is the extraction of the required relative phases from the given interference
    measurements.  Writing
    $\fourier x (\omega) = \absn{\fourier x (\omega)} \, \e^{\I \phi (\omega)}$, we can
    rearrange the first interference measurement to
    \begin{align*}
      &\abs{\fourier x \mleft( \omega \mright)
        + \e^{\I \alpha_1} \, \fourier x \mleft( \omega - \mu \mright)}^2
      \\[\fskip]
      & \qquad= \abs{\fourier x \mleft( \omega \mright)}^2
        + \abs{\fourier x \mleft( \omega - \mu \mright)}^2
        + 2 \abs{\fourier x \mleft( \omega \mright)}
        \abs{\fourier x \mleft( \omega - \mu \mright)}
        \Re \mleft[ \e^{\I ( \phi ( \omega - \mu ) - \phi ( \omega ) + \alpha_1 )} \mright].
    \end{align*}
    Replacing $\alpha_1$ by $\alpha_2$, we obtain a similar representation for the second
    interference measurement.  Consequently, if both moduli $\absn{\fourier x (\omega)}$
    and $\absn{\fourier x (\omega - \mu)}$ of the Fourier transform $\fourier x$ are
    non-zero, we can determine
    \begin{equation*}
        \Re \mleft[ \e^{\I ( \phi(\omega - \mu) - \phi(\omega) + \alpha_1)} \mright]
        \quad \text{and} \quad
        \Re \mleft[ \e^{\I ( \phi(\omega - \mu) - \phi(\omega) + \alpha_2)} \mright].
    \end{equation*}

    In order to extract the phase difference $\phi(\omega-\mu) - \phi(\omega)$ from these
    two values, we apply Euler's formula and afterwards the addition theorem for cosine.
    In this manner, we obtain
    \begin{align*}
      \Re \mleft[ \e^{\I (\phi(\omega - \mu) - \phi(\omega) + \alpha_1)} \mright]
      &= \cos\mleft( \alpha_1 \mright)
        \cos\mleft( \phi \mleft( \omega - \mu \mright)
        - \phi \mleft( \omega \mright) \mright)
      \\[\fsmallskip]
      & \qquad-
        \sin\mleft( \alpha_1 \mright)
        \sin\mleft( \phi \mleft( \omega - \mu \mright)
        - \phi \mleft( \omega \mright) \mright)
        \shortintertext{and}
        \Re \mleft[ \e^{\I (\phi(\omega - \mu) - \phi(\omega) + \alpha_2)} \mright]
      &= \cos\mleft( \alpha_2 \mright)
        \cos\mleft( \phi \mleft( \omega - \mu \mright)
        - \phi \mleft( \omega \mright) \mright)
      \\[\fsmallskip]
      & \qquad-
        \sin\mleft( \alpha_2 \mright)
        \sin\mleft( \phi \mleft( \omega - \mu \mright)
        - \phi \mleft( \omega \mright) \mright).
    \end{align*}
    Since the values on the left-hand side are known, we can consequently determine the
    sine and cosine of the wanted phase difference $\phi(\omega-\mu) - \phi(\omega)$ by
    solving a simple linear equation system.  Here, the determinant of the occurring
    matrix is given by
    \begin{equation*}
        \det
        \begin{psmallmatrix}
            \cos \alpha_1 & -\sin \alpha_1 \\
            \cos \alpha_2 & -\sin \alpha_2
        \end{psmallmatrix}
        = \sin \alpha_1 \, \cos \alpha_2 - \cos \alpha_1 \, \sin \alpha_2
        = \sin \mleft( \alpha_1 - \alpha_2 \mright),
    \end{equation*}
    which ensures a unique solution whenever $\alpha_1-\alpha_2$ does not coincide with a
    multiple of $\pi$.  Consequently, we can always determine the required phase
    difference $\phi(\omega-\mu)-\phi(\omega)$ for a certain $\omega$.  With the extracted
    phase difference, the theorem can now be justified as discussed in the proof of
    Theorem~\ref{the:interf-polar-id}.  \qed%
\end{Proof}

\begin{Remark}
    As mentioned before, Cand\`es \etal\ consider a slightly different phase retrieval
    problem in \cite[Theorem~3.1]{CESV13}.  More precisely, they deal with the problem to
    recover a finite-dimensional vector from the intensities of its discrete Fourier
    transform.  Using our notation, we can state this problem as follows: recover the
    signal $x$ whose support of length $N$ is contained in the interval from $0$ to $M-1$
    from the (discrete-time) Fourier intensities
    \begin{equation*}
        \abs{\fourier x \mleft( \tfrac{2 \pi k}{M} \mright)} \!,
        \quad
        \abs{\fourier x \mleft( \tfrac{2 \pi k}{M} \mright)
            + \fourier x \mleft( \tfrac{2 \pi \left(k-\ell\right)}{M} \mright)} \!,
        \quad \text{and} \quad
        \abs{\fourier x \mleft( \tfrac{2 \pi k}{M} \mright)
            - \I \, \fourier x \mleft( \tfrac{2 \pi \left(k-\ell\right)}{M} \mright)}
    \end{equation*}
    for all integers $k=0,\dots,M-1$ and for a certain integer $\ell$.  Under the
    additional assumption that $\ell$ and $M$ are co-prime, and that the given Fourier
    samples $\absn{\fourier x (\nicefrac{2\pi k}{M})}$ are non-zero, Cand\`es \etal\
    show that the unknown signal $x$ can be uniquely recovered.

    Recalling that the Fourier intensity of a discrete-time signal with support length $N$
    is entirely determined by $2N-1$ arbitrary samples, we can directly compare
    Theorem~\ref{the:two-interf-arb-rot} for $\alpha_1 \coloneqq 0$,
    $\alpha_2 \coloneqq -\nicefrac{\pi}{2}$, and
    $\mu \coloneqq 2 \pi \, \nicefrac{\ell}{M}$ with the results of Cand\`es \etal\ and
    see that both statements are almost identical.  In the special case that $N$ and $M$
    coincide, the main difference between both statements is that we need twice as many
    measurements. Anyway, this enables us to neglect the assumption that the given samples
    of the Fourier intensity have to be non-zero.  Further, we can determine the unknown
    position of the current support completely from the given Fourier intensities.

    Finally, the integers $N$ and $M$ have a slightly different meaning.  With the
    dimension $M$, we determine the interval $\{0,\dots,M-1\}$ that contains the non-zero
    components of the considered signal $x$.  The current support length $N$ of this
    signal can however be much smaller than the assumed dimension $M$.  Consequently, if
    $M$ is only a rough estimation, then Theorem~\ref{the:two-interf-arb-rot} allows us to
    recover the wanted signal from a much smaller set of measurements.  \qed%
\end{Remark}

\section{Continuous-time phase retrieval}
\label{sec:cont-time-phase}

\subsection{Characterization of the occurring ambiguities}
\label{sec:phase-retr-arbitr-sig}

Different from the previous section, we now consider the continuous-time phase retrieval
problem, which has a completely different behavior as the discrete-time equivalent.  Here
we are faced with the recovery of a continuous-time signal or a function
$f \colon \R \rightarrow \C$ in $L^2$ with compact support from its Fourier intensity
$\absn{\Fourier [f]}$, where the corresponding Fourier transform is given by
\begin{equation*}
    \fourier f \mleft( \omega \mright)
    \coloneqq \Fourier \mleft[ f \mright] \mleft( \omega \mright)
    \coloneqq \int_{-\infty}^\infty f \mleft( t \mright) \, \e^{-\I \omega t} \diff t.
\end{equation*}

The ambiguousness of this problem has been studied by Akutowicz \cite{Aku56,Aku57},
Walther \cite{Wal63}, and Hofstetter \cite{Hof64}, for instance.  To reduce the set of
possible ambiguities, we will again employ interference measurements with different kinds
of reference signals.  For this purpose, we require a suitable characterization of the
solution set.  Here we follow the approach by Hofstetter \cite{Hof64}, where the solutions
are presented as an infinite product with respect to the zero set of the analytic
continuations of their Fourier transforms.

In our case, the continuation of $\Fourier [f]$ is given through the theorem of
Paley-Wiener, see for instance \cite[ Theorem~6.8.1]{Boa54}.  Using the (two-sided)
Laplace transform, we can write this well-known theorem in the following form.

\begin{Theorem}[Paley-Wiener]
    \label{the:Paley-Wiener}%
    The function $F \colon \C \rightarrow \C$ is an entire function of exponential type
    and belongs to $L^2$ on the imaginary axis if and only if $F$ is the Laplace transform
    \begin{equation*}
        F \mleft(\zeta\mright)
        \coloneqq
        \Laplace \mleft[ f \mright]\mleft( \zeta \mright)
        \coloneqq
        \int_{-\infty}^\infty f\mleft( t \mright) \, \e^{-\zeta t} \diff t,
    \end{equation*}
    of a function $f$ in $L^2$ with compact support.%
\end{Theorem}

\vspace*{-6pt}

In this context, a function $F \colon \C \rightarrow \C$ is called \emph{entire} if the
function is analytic over the whole complex plane.  If $F$ moreover grows no faster than
an exponential, which means that $F$ can be bounded by
\vspace*{-4pt}
\begin{equation*}
    \abs{F \mleft( \zeta \mright)}
    \le A \, \e^{B \abs{\zeta}},
    \addmathskip
    \vspace*{-4pt}
\end{equation*}
then the entire function $F$ is of \emph{exponential type}, see for instance
\cite[\p~53]{You80}.  As a consequence of the theorem of Paley-Wiener, the Laplace
transform $F$ is now the unique analytic continuation of the Fourier transform
$\Fourier [f]$ from the imaginary axis to the complex plane because the restriction
$F(\I \cdot)$ of the holomorphic Laplace transform $F$ to the imaginary axis obviously
coincides with $\Fourier [f]$.

The next key instruments in the characterization of the occurring ambiguities are the
\emph{autocorrelation signal}
\vspace*{-4pt}
\begin{equation*}
    a \mleft( t \mright)
    \coloneqq \int_{-\infty}^\infty \overline{f \mleft( s \mright)} \,
    f \mleft( s + t \mright) \diff s
    \addmathskip
    \vspace*{-4pt}
\end{equation*}
and the \emph{autocorrelation function} $A \coloneqq \Laplace [ a ]$.  Since the
autocorrelation signal of a compactly supported function $f$ in $L^2$ is again a compactly
supported function in $L^2$, the autocorrelation function $A$ can be interpreted as the
analytic continuation of the Fourier transform $\fourier a$.  Further, the autocorrelation
is closely related to the given Fourier intensity, see \cite{Hof64}.

\vspace*{-6pt}

\begin{Proposition}
    \label{prop:four-int-cont-autocor}%
    Let $f$ be a continuous-time signal in $L^2$ with compact support. Then the
    autocorrelation function $A$ is the analytic continuation of the squared Fourier
    intensity $\absn{\Fourier[f]}^2$ from the imaginary axis to the complex plane.%
\end{Proposition}

\vspace*{-10pt}

\begin{Proof}
    Using the definition of the autocorrelation signal, we can write the autocorrelation
    function $A$ of the considered signal $f$ as
    \vspace*{-5pt}
    \begin{equation}
        \begin{aligned}
            A \mleft( \zeta \mright)
            &= \int_{-\infty}^\infty \int_{-\infty}^\infty
            \overline{f \mleft( s \mright)} \, f \mleft( s + t \mright) \, \e^{-\zeta t}
            \diff s \diff t
            \\[\fskip]
            &= \int_{-\infty}^\infty
            \int_{-\infty}^\infty \overline{f \mleft( s \mright)} \, f \mleft( t \mright)
            \, \e^{-\zeta t} \, \e^{\zeta s} \diff t \diff s = F \bigl( \zeta \bigr) \,
            \overline{F \bigl( -\overline \zeta \bigr)}.
        \end{aligned}
        \label{eq:four-int-cont-autocor}
        \vspace*{0.1pt}
    \end{equation}
    If we now consider the restriction of this identity to the imaginary axis, we have
    \begin{equation*}
        A \mleft( \I \omega \mright)
        = F \mleft( \I \omega \mright) \, \overline{ F \mleft( \I\omega \mright) }
        = \fourier f \mleft( \omega \mright) \, \overline{ \fourier f \mleft( \omega
            \mright) }
        = \absm{\fourier f \mleft( \omega \mright)}^2 \!,
        \addmathskip
    \end{equation*}
    which implies that the restriction of the autocorrelation function $A$ coincides with
    the squared Fourier intensity $\absn{\Fourier [f]}^2$ as claimed.  Since $A$ is a
    holomorphic function by the theorem of Paley-Wiener, the assertion follows.  \qed%
\end{Proof}

After this preliminaries, we recall the characterization of all occurring ambiguities in
the continuous-time phase retrieval problem, \cite[Theorem~I]{Hof64}.

\begin{Theorem}[Hofstetter]
    \label{the:char-amb-Hofstetter}%
    Let $f$ be a continuous-time signal in $L^2$ with compact support.  Then the
    Laplace transform of each continuous-time signal $g$ in $L^2$ with compact support
    and the same Fourier intensity $\absn{\Fourier [g]} = \absn{\Fourier [f]}$ can be
    written in the form
    \begin{equation*}
        G \mleft( \zeta \mright)
        = C \, \zeta^m \e^{\zeta \gamma} \prod_{j=1}^\infty
        \left( 1 - \tfrac{\zeta}{\eta_j} \right) \e^{\frac{\zeta}{\eta_j}}
        \addmathskip
    \end{equation*}
    where the absolute value $\absn{C}$ and the imaginary part $\Im \gamma$ of the complex
    constants $C$ and $\gamma$ coincide for all signals $g$, and where $\eta_j$ is chosen
    from the zero pair $(\xi_j, -\overline \xi_j)$ of the autocorrelation function $A$.%
\end{Theorem}

\begin{Proof}
    For the sake of convenience for the reader, we give a short proof following the lines
    of \cite{Hof64}.  Let $f$ and $g$ be signals with compact support in $L^2$, and let
    $F$ and $G$ be their Laplace transforms being entire functions of exponential type.
    By Hadamard's factorization theorem \cite{Boa54}, we can represent the entire
    functions $F$ and $G$ by
    \begin{equation*}
      F \mleft( \zeta \mright)
      =  C_1 \, \zeta^{m_1} \, \e^{\zeta \gamma_1 }
        \prod_{j=1}^{\infty} \left( 1 - \tfrac{\zeta}{\xi_j} \right)
        \e^{\frac{\zeta}{\xi_j}}
        \quad\text{and}\quad
        G \mleft( \zeta \mright)
      =  C_2 \, \zeta^{m_2} \, \e^{\zeta \gamma_2}
        \prod_{j=1}^{\infty} \left( 1 - \tfrac{\zeta}{\eta_j} \right)
        \e^{\frac{\zeta}{\eta_j}}
    \end{equation*}
    with respect to their non-vanishing zeros $\xi_j$ and $\eta_j$.

    Assuming that $\absn{\Fourier [g]} = \absn{\Fourier [f]}$, we can conclude that the
    autocorrelation functions of $f$ and $g$ coincide.  Using
    \eqref{eq:four-int-cont-autocor}, we can therefore represent the common
    autocorrelation function $A$ in terms of $F$ or $G$ by
    \begin{equation}
        \label{eq:fact-auto-F-G}
        A \mleft( \zeta \mright)
        = F \bigl( \zeta \bigr) \, \overline{F \bigl( - \overline \zeta \bigr)}
        = G \bigl( \zeta \bigr) \, \overline{G \bigl( - \overline \zeta \bigr)}.
    \end{equation}
    Due to this factorization, all zeros of the autocorrelation function $A$ obviously
    occur in pairs of the form \raisebox{0pt}[0pt][0pt]{$(\xi_j^{\,}, -\overline \xi_j)$},
    where $\xi_j$ is a zero of $F$.  Since an analogous observation follows from the
    factorization of $G$, we can resort the zeros $\eta_j$ so that
    $\eta_j^{\,} = \xi_j^{\,}$ or
    \raisebox{0pt}[0pt][0pt]{$\eta_j^{\,} = - \overline \xi_j$}.  Further, the
    multiplicities $m_1$ and $m_2$ of the zero at the origin must also be equal.
    Consequently, the possibly infinite products in the factorizations of $A$ coincide,
    and we can reduce \eqref{eq:fact-auto-F-G} to
    \begin{equation*}
        \abs{C_1}^2 \e^{\zeta (\gamma_1^{\,} - \overline \gamma_1)}
        = \abs{C_2}^2 \e^{\zeta (\gamma_2^{\,} - \overline \gamma_2)}
        \quad\text{or}\quad
        \abs{C_1}^2 \e^{2 \zeta \Im[\gamma_1] }
        = \abs{C_2}^2 \e^{2 \zeta \Im[\gamma_2] },
        \addmathskip
    \end{equation*}
    which shows that the absolute values $\absn{C_1}$ and $\absn{C_2}$ and also the
    imaginary parts $\Im \gamma_1$ and $\Im \gamma_2$ coincide.  \qed%
\end{Proof}

Similarly to the discrete-time setting, we divide the occurring ambiguities characterized
by Theorem~\ref{the:char-amb-Hofstetter} in two different classes.  Since the rotation, the time
shift by a real number, and the reflection and conjugation of a solution always result in
a further solution of the considered phase retrieval problem, \cf\
Proposition~\ref{prop:trivial-amb}, we call these ambiguities trivial.  The remaining ambiguities
are non-trivial.  Using the characterization in Theorem~\ref{the:char-amb-Hofstetter}, we can
also represent the ambiguities in the following manner.

\begin{Proposition}
    \label{prop:repr-cont-sol-conv}%
    Let $f$ and $g$ be two continuous-time signals in $L^2$ with compact support and the
    same Fourier intensity $\absn{\Fourier [f]} = \absn{\Fourier [g]}$.  Then there exist two entire
    functions $F_1$ and $F_2$ of exponential type such that
    \begin{equation*}
        F \mleft( \zeta \mright)
        = F_1 \mleft( \zeta \mright) \, F_2 \mleft( \zeta \mright)
        %\submathskip
    \end{equation*}
    and
    \begin{equation*}
        G \mleft( \zeta \mright)
        = \e^{\I \alpha} \, \e^{-\zeta t_0} \, \overline{ F_1 \bigl(-\overline \zeta
            \bigr)} \, F_2  \bigl( \zeta \bigr),
        \addmathskip
    \end{equation*}
    where $\alpha$ and $t_0$ are suitable real numbers.%
\end{Proposition}

\begin{Proof}
    Applying Theorem~\ref{the:char-amb-Hofstetter}, we can represent the Laplace transforms $F$
    and $G$ by
    \begin{equation*}
        F \mleft( \zeta \mright)
        = C_1 \, \zeta^m \, \e^{\zeta \gamma_1} \prod_{j=1}^\infty
        \left( 1 - \tfrac{\zeta}{\xi_j} \right) \e^{\frac{\zeta}{\xi_j}}
        \quad\text{and}\quad
        G \mleft( \zeta \mright)
        = C_2 \, \zeta^m \, \e^{\zeta \gamma_2} \prod_{j=1}^\infty
        \left( 1 - \tfrac{\zeta}{\eta_j} \right) \e^{\frac{\zeta}{\eta_j}}
    \end{equation*}
    with $\absn{C_1}=\absn{C_2}$, $\Im \gamma_1 = \Im \gamma_2$, and
    $\eta_j \in (\xi_j, -\overline \xi_j)$.  Now, we can resort the zeros of $G$ such that
    \begin{equation*}
        \eta_j =
        \begin{cases}
            -\overline \xi_j & \xi_j^{\,} \in \Lambda, \\
            \phantom{-}\xi_j^{\,} & \text{else}
        \end{cases}
        \addmathskip
    \end{equation*}
    for an appropriate subset $\Lambda$ of the zero set
    $\Xi\coloneqq \{ \xi_j \colon j \in \N \}$ of $F$.

    Based on this subset, we define the two possibly infinite products $F_1$ and $F_2$ by
    \begin{equation*}
        F_1 \mleft( \zeta \mright)
        \coloneqq \prod_{\xi_j \in \Lambda} \left( 1 - \tfrac{\zeta}{\xi_j} \right)
        \e^{\frac{\zeta}{\xi_j}}
        \quad\text{and}\quad
        F_2 \mleft( \zeta \mright)
        \coloneqq C_1 \, \zeta^m \, \e^{\zeta \gamma_1}
        \smashoperator{\prod_{\xi_j \in \Xi \setminus\Lambda}} \left( 1 - \tfrac{\zeta}{\xi_j} \right)
        \e^{\frac{\zeta}{\xi_j}}.
    \end{equation*}
    Due to the fact that the convergence exponent -- the infimum of the positive numbers
    $\alpha$ for which the series $\sum_{j=1}^\infty \absn{\xi_j}^{-\alpha}$ converges --
    of the zeros of an entire function is always less than or equal to the order of the
    entire function, see for instance \cite[Theorem~2.5.18]{Boa54}, the zeros $\xi_j$ can
    at most have the convergence exponent one.  Since the zeros of $F_1$ and $F_2$ are
    merely subsets of the zeros of $F$, the corresponding convergence exponents thus also
    have to be less than or equal to one.  Borel's theorem now implies that the possibly
    infinite products $F_1$ and $F_2$ are again entire functions of exponential type, see
    \cite{Mar77}.

    Obviously, we have the factorization $F = F_1 \, F_2$.  In order to achieve the
    factorization of $G$, we consider the reflection of the first factor given by
    \begin{equation*}
        \overline{F_1 \mleft( - \overline \zeta \mright)}
        = \prod_{\xi_j \in \Lambda} \, \Bigl( 1 - \tfrac{\zeta}{-\overline\xi_j} \Bigr)
        \, \e^{\frac{\zeta}{-\overline\xi_j}}.
    \end{equation*}
    Hence, the reflection $\overline{F_1(- \bar \cdot)}$ possesses the zeros
    $\eta_j = -\overline \xi_j$ for all $\xi_j$ in $\Lambda$, which implies that the zeros
    of the product \raisebox{0pt}[0pt][0pt]{$\overline{F_1(- \bar \cdot)} \, F_2$} and $G$
    coincide.

    Finally, since the absolute values $\absn{C_1}$ and $\absn{C_2}$ and the imaginary
    parts $\Im \gamma_1$ and $\Im \gamma_2$ have to be equal by
    Theorem~\ref{the:char-amb-Hofstetter}, the entire functions
    \raisebox{0pt}[0pt][0pt]{$\overline{F_1(- \bar \cdot)} \, F_2$} and $G$ can only
    differ by a rotation $\e^{\I\alpha}$ and a time shift $\e^{-\zeta t_0}$.  Choosing the
    real numbers $\alpha$ and $t_0$ suitably, we obtain the wanted factorization in the
    assertion.  \qed%
\end{Proof}

\begin{Remark}
    \label{rem:repr-cont-sol-conv:2}
    We can interpret Proposition~\ref{prop:repr-cont-sol-conv} in a slightly different way: if the
    restrictions of the entire functions $F$ and $G$ of exponential type to the imaginary
    axis coincide, then we can always find two entire functions $F_1$ and $F_2$ to
    factorize $F$ and $G$ in the manner of Proposition~\ref{prop:repr-cont-sol-conv}.  This
    observation can now be generalized to entire function of arbitrary order, see
    \cite[Lemma~1]{Mar14}.  \qed%
\end{Remark}

\begin{Remark}
    For the discrete-time setting, an equivalent statement of
    Proposition~\ref{prop:repr-cont-sol-conv} holds where the Laplace transforms $F$ and $G$ are
    replaced by the discrete-time Fourier transforms of the two signals.  Moreover, for
    two discrete-time signals $x$ and $y$ with $\absn{\fourier x} = \absn{\fourier y}$,
    there always exists two discrete-time signals $x_1$ and $x_2$ with finite support such
    that
    \begin{equation*}
        x = x_1 * x_2 \coloneqq \Bigl( \sum_{k\in\Z} x_1[n-k] \, x_2[k] \Bigr)_{n\in\Z}
    \qquad\text{and}\qquad
        y = \e^{\I \alpha} \, x_1 [\cdot - n_0] * \overline{x_2 [- \cdot]}
    \end{equation*}
    with suitable $\alpha \in \R$ and $n_0 \in \Z$, see \cite[Theorem~2.3]{BP15}.  In
    other words, each ambiguity of the discrete-time phase retrieval problem can be
    represented by an appropriate convolution and a suitable rotation, shift, and
    reflection and conjugation of the occurring factors.  \qed%
\end{Remark}

\subsection{Ensuring uniqueness in the continuous-time phase retrieval}
\label{sec:ensur-uniq-cont}

Differently from the discrete-time phase retrieval problem, the continuous-time version to
recover a certain signal usually possesses infinitely many non-trivial ambiguities.
Hence, we are once more faced with the question: how can we ensure the unique recovery of
the unknown signal, or how can we at least reduce the occurring ambiguities to a
sufficiently small set.

Like in the discrete-time setting, we will employ different kinds of interference
measurements.  Here Klibanov \etal\ show that the additional interference measurement with
a known reference can almost ensure the unique recovery of a distribution with compact
support.  Adapting this observation \cite[Proposition~6.5]{KST95} to the continuous-time
phase retrieval problem, we obtain the following statement.

\begin{Proposition}[Klibanov \etal]
    \label{prop:cont-known-ref}%
    Let $f$ and $h$ be two continuous-time signals in $L^2$ with compact support, where
    the non-vanishing reference signal $h$ is known beforehand.  Then the signal $f$ can
    be recovered from the Fourier intensities
    \begin{equation*}
        \abs{\Fourier[f]}
        \quad\text{and}\quad
        \abs{\Fourier [f+h]}
    \end{equation*}
    except for at most one ambiguity.%
\end{Proposition}

\begin{Proof}
    We give a new proof of the statement by adapting the proof of the discrete-time
    equivalent in \cite[Theorem~4.3]{BP15}.  Writing the Fourier transforms of the signals
    $f$ and $h$ in their polar representations
    \begin{equation*}
        \Fourier [f] \mleft(\omega\mright) = \abs{\Fourier [f] \mleft(\omega\mright)}
        \e^{\I \phi (\omega)}
        \quad\text{and}\quad
        \Fourier [h] \mleft(\omega\mright) = \abs{\Fourier [h] \mleft(\omega\mright)}
        \e^{\I \psi (\omega)},
    \end{equation*}
    where $\phi$ and $\psi$ denote the corresponding phase functions, we can represent the
    given interference measurement by
    \begin{equation*}
        \absn{\Fourier [f+h](\omega)}^2
        = \absn{\Fourier [f](\omega)}^2 + \absn{\Fourier [h](\omega)}^2 + 2 \absn{\Fourier
            [f](\omega)} \absn{\Fourier [h](\omega)} \cos ( \phi(\omega) - \psi(\omega)).
    \end{equation*}
    In other words, the phase difference $\phi-\psi$ is given by
    \begin{equation}
        \label{eq:cont-known-ref:phase-diff}
        \phi \mleft( \omega \mright) - \psi \mleft( \omega \mright)
        = \pm \arccos \left( \frac{\absn{\Fourier \mleft[ f+h \mright] \mleft( \omega
                    \mright)}^2 - \absn{\Fourier \mleft[ f \mright] \mleft( \omega
                    \mright)}^2 - \absn{\Fourier \mleft[ h \mright] \mleft( \omega
                    \mright)}^2}{2 \absn{\Fourier \mleft[ f \mright] \mleft( \omega
                    \mright)} \, \absn{\Fourier \mleft[ h \mright] \mleft( \omega \mright)}}
        \right) + 2 \pi k
    \end{equation}
    with an appropriate integer $k$ whenever $\Fourier [f](\omega)$ and
    $\Fourier [h](\omega)$ are non-zero.

    Since $\Fourier [f]$ and $\Fourier [h]$ are continuous, we can find a small interval
    for $\omega$ where the sign in \eqref{eq:cont-known-ref:phase-diff} does not change
    and the integer $k$ is fixed.  Further, since $\phi-\psi$ is the phase function of the
    product \raisebox{0pt}[0pt][0pt]{$\Fourier [f] \overline{\Fourier [h]}$}, which is the
    restriction of the entire function
    \raisebox{0pt}[0pt][0pt]{$F( \cdot) \, \overline{H ( - \bar \cdot)}$} to the imaginary
    axis, we can extend $\phi-\psi$ uniquely from the interval to the complete frequency
    domain. Consequently, there exist at most two distinct phase differences $\phi-\psi$.
    Here the phase $\breve \phi$ of the second solution $\breve f$ and the phase $\phi$
    are related by
    \begin{equation*}
        \phi \mleft( \omega \mright) - \psi \mleft( \omega \mright)
        = - \breve \phi \mleft( \omega \mright) + \psi \mleft( \omega \mright) + 2\pi k.
        \addmathskip
    \end{equation*}
    Hence, the Fourier transform of the second solution has to be of the form
    \begin{equation}
        \label{eq:cont-known-ref:sec-sol}
        \Fourier [\breve f] = \abs{\Fourier [f]} \e^{-\I\phi + 2 \I \psi},
        \submathskip
    \end{equation}
    which completes the proof.  \qed%
\end{Proof}

\begin{Remark}
    The main benefit of the proof of Proposition~\ref{prop:cont-known-ref} given above is that we
    obtain an explicit representation of the second possible solution in dependence of the
    phase of the reference signal.  Considering the Fourier transform
    \eqref{eq:cont-known-ref:sec-sol}, we can have doubts whether the corresponding
    continuous-time signal $\breve f$ is really a signal with compact support and hence a
    valid solution of the problem.  Indeed, the Fourier transform
    \eqref{eq:cont-known-ref:sec-sol} does not have to be the restriction of an entire
    function and does not even have to be a continuous function at all because the phase
    $\psi$ of the continuous function $\Fourier [h]$ itself can possesses discontinuities.
    Further, the theorem of Paley-Wiener here implies that the second solution $\breve f$
    does not have to have a compact support and may thus be an invalid solution of the
    considered continuous-time phase retrieval problem.  \qed%
\end{Remark}

Next, we replace the known reference signal $h$ within the interference $f+h$ by an
unknown reference.  Inspired by the discrete-time equivalent in \cite[Theorem~4.4]{BP15},
we show that the continuous-time signal $f$ and the unknown reference $h$ are uniquely
determined by the Fourier intensities of $f$, $h$, and $f+h$ up to common trivial
ambiguities.  This means that we can recover $f$ and $h$ up to common rotations or time
shifts or up to the reflection and conjugation of the two signals.

\begin{Theorem}
    \label{the:cont-unknown-ref}%
    Let $f$ and $h$ be two continuous-time signals in $L^2$ with compact support.  If the
    non-vanishing zeros of the Laplace transformed signal $F$ and $H$ form disjoint sets,
    then both signals $f$ and $h$ can be recovered from the Fourier intensities
    \begin{equation*}
        \absn{\Fourier \mleft[ f \mright]},
        \quad
        \absn{\Fourier \mleft[ h \mright]},
        \quad\text{and,}\quad
        \absn{\Fourier \mleft[ f+h \mright]}
    \end{equation*}
    uniquely up to common trivial ambiguities.%
\end{Theorem}

\begin{Proof}
    Let $\breve f$ and $\breve h$ be a further solution pair of the considered problem
    with
    \begin{equation*}
        \absn{\Fourier [f]} = \absn{\Fourier [\breve f]},
        \quad
        \absn{\Fourier [h]} = \absn{\Fourier [\breve h]},
        \quad\text{and}\quad
        \absn{\Fourier [f+h]} = \absn{\Fourier [\breve f + \breve h]}.
    \end{equation*}
    Applying Proposition~\ref{prop:repr-cont-sol-conv}, we can represent the two solution pairs in
    the frequency domain by an appropriate factorization of the Laplace transform $F$
    and $H$ of the original signals.  In this manner, we obtain the factorizations
    \begin{align*}
      F \mleft( \zeta \mright)
      = F_1 \mleft( \zeta \mright) \, F_2 \mleft( \zeta \mright)
      &\quad\text{and}\quad
        \breve F \mleft( \zeta \mright)
        = \e^{\I\alpha_1} \, \e^{-\zeta t_1} \,
        \overline{F_1 \bigl( - \overline \zeta \bigr)} \, F_2 \bigl( \zeta \bigr)
        \intertext{and further}
        H \mleft( \zeta \mright)
        = H_1 \mleft( \zeta \mright) \, H_2 \mleft( \zeta \mright)
      &\quad\text{and}\quad
        \breve H \mleft( \zeta \mright)
        = \e^{\I\alpha_2} \, \e^{-\zeta t_2} \,
        \overline{H_1 \bigl( - \overline \zeta \bigr)} \, H_2 \bigl( \zeta \bigr)
    \end{align*}
    for some real numbers $\alpha_1$, $\alpha_2$, $t_1$, $t_2$ and entire functions $F_1$,
    $F_2$, $H_1$, $H_2$.

    In the next step, we consider the analytic continuation of the squared interference
    measurement or the corresponding autocorrelation function, see
    Proposition~\ref{prop:four-int-cont-autocor}.  With the representation in
    \eqref{eq:four-int-cont-autocor}, we can write the given interference measurement as
    \begin{equation*}
        \left( F \bigl( \zeta \bigr) + H \bigl( \zeta \bigr) \right)
        \left( \overline{F \bigl( - \overline \zeta \bigr)} + \overline{H \bigl( -
                \overline \zeta \bigr)} \right)
        = \left( \breve F \bigl( \zeta \bigr) + \breve H \bigl( \zeta \bigr) \right) 
        \left( \overline{ \breve F \bigl( - \overline \zeta \bigr)} + \overline{ \breve H
                \bigl( - \overline \zeta \bigr)} \right) 
    \end{equation*}
    or in the simplified form
    \begin{equation*}
        F \bigl( \zeta \bigr) \, \overline{H \bigl( - \overline \zeta \bigr)}
        + \overline{ F \bigl( - \overline \zeta \bigr)} \, H \bigl( \zeta \bigr)
        = \breve F \bigl( \zeta \bigr) \, \overline{\breve H \bigl( - \overline \zeta
            \bigr)} + \overline{ \breve F \bigl( - \overline \zeta \bigr)} \, \breve H
        \bigl( \zeta \bigr).
    \end{equation*}
    Incorporating the found factorizations of $F$ and $H$, we obtain
    \begin{align*}
      & F_1 \bigl( \zeta \bigr) \, F_2 \bigl( \zeta \bigr) \, \overline{H_1 \bigl(-\overline
        \zeta \bigr)} \, \overline{H_2 \bigl(- \overline \zeta \bigr)}
        + \overline{F_1 \bigl(-\overline \zeta\bigr)} \, \overline{F_2 \bigl( - \overline
        \zeta \bigr)} \, H_1 \bigl( \zeta \bigr) \, H_2 \bigl( \zeta \bigr)
      \\[\fskip]
      & \quad = \e^{\I(\alpha_1 - \alpha_2)} \, \e^{-\zeta(t_1 - t_2)} \overline{F_1
        \bigl(-\overline \zeta \bigr)} \, F_2 \bigl( \zeta \bigr) \, H_1 \bigl( \zeta
        \bigr) \, \overline{H_2 \bigl( - \overline \zeta \bigr)}
      \\[\fsmallskip]
      & \quad\quad
        + \e^{\I(\alpha_2 - \alpha_1)} \, \e^{-\zeta(t_2 - t_1)} F_1 \bigl( \zeta \bigr)
        \, \overline{F_2 \bigl(- \overline \zeta \bigr)} \, \overline{ H_1 \bigl( -
        \overline \zeta \bigr)} \, H_2 \bigl( \zeta \bigr)
    \end{align*}
    and thus
    \begin{equation}
        \label{eq:cont-unknown-ref:fact-intf}
        \begin{aligned}
            & \left[ \e^{-\I\alpha_1} \, \e^{\zeta t_1} \, F_1 \bigl( \zeta \bigr) \,
                \overline{H_1 \bigl(-\overline \zeta\bigr)} - \e^{-\I\alpha_2} \,
                \e^{\zeta t_2} \, \overline{F_1 \bigl(-\overline \zeta \bigr)} \, H_1
                \bigl( \zeta \bigr) \right]
            \\[\fsmallskip]
            &\quad \cdot \left[ \e^{\I\alpha_1} \, \e^{- \zeta t_1} \, F_2 \bigr( \zeta
                \bigr) \, \overline{H_2 \bigl(- \overline \zeta \bigr)} - \e^{\I \alpha_2}
                \, \e^{- \zeta t_2} \, \overline{F_2\bigl( - \overline \zeta \bigr)} \,
                H_2 \bigl( \zeta \bigr) \right] = 0.
        \end{aligned}
    \end{equation}

    Remembering that $F_1$, $F_2$, $H_1$, and $H_2$ are entire functions, we observe that
    both factors in \eqref{eq:cont-unknown-ref:fact-intf} are entire functions too, and
    that at least one of both factors thus has to be constantly zero.  In order to
    investigate the two different cases more precisely, we employ the explicit
    construction of the entire functions $F_1$ and $F_2$ in the proof of
    Proposition~\ref{prop:repr-cont-sol-conv}.  Using a similar procedure for $H_1$ and $H_2$, and
    denoting the sets of all non-vanishing zeros of $F$ and $H$ by $\Xi_1$ and $\Xi_2$
    respectively, we can represent the four functions by
    \begin{align*}
      F_1 \mleft( \zeta \mright)
      = \prod_{\xi_j \in \Lambda_1} \left( 1 - \tfrac{\zeta}{\xi_j} \right)
      \e^{\frac{\zeta}{\xi_j}}
      &\quad\text{and}\quad
        F_2 \mleft( \zeta \mright)
        = C_1 \, \zeta^{m_1} \, \e^{\zeta \gamma_1}
        \smashoperator{\prod_{\xi_j \in \Xi_1\setminus\Lambda_1}}
        \left( 1 - \tfrac{\zeta}{\xi_j} \right) \e^{\tfrac{\zeta}{\xi_j}}
        \intertext{and further}
        H_1 \mleft( \zeta \mright)
        = \prod_{\eta_j \in \Lambda_2} \left( 1 - \tfrac{\zeta}{\eta_j} \right)
        \e^{\frac{\zeta}{\eta_j}}
      &\quad\text{and}\quad
        H_2 \mleft( \zeta \mright)
        = C_2 \, \zeta^{m_2} \, \e^{\zeta \gamma_2}
        \smashoperator{\prod_{\eta_j \in \Xi_2\setminus\Lambda_2}}
        \left( 1 - \tfrac{\zeta}{\eta_j} \right) \e^{\tfrac{\zeta}{\eta_j}} \!,
    \end{align*}
    where $\Lambda_1$ and $\Lambda_2$ are appropriate subsets of $\Xi_1$ and $\Xi_2$.

    In the following, we firstly assume that the second factor of
    \eqref{eq:cont-unknown-ref:fact-intf} is zero, which directly implies that the
    equation
    \begin{equation}
        \label{eq:cont-unknown-ref:sec-fact}\hspace*{-15pt}
        \begin{aligned}
            &\left(-1\right)^{m_2} C_1 \, \overline C_2 \, \e^{\I\alpha_1} \, 
            \zeta^{m_1+m_2} \, \e^{-\zeta (t_1-\gamma_1+\overline \gamma_2)}
            \smashoperator{\prod_{\xi_j \in \Xi_1\setminus\Lambda_1}}
            \, \Bigl( 1 - \tfrac{\zeta}{\xi_j} \Bigr) \,
            \e^{\frac{\zeta}{\xi_j}}
            \smashoperator{\prod_{\eta_j \in \Xi_2\setminus\Lambda_2}} \,
            \Bigl( 1 - \tfrac{\zeta}{-\overline\eta_j} \Bigr) \,
            \e^{\frac{\zeta}{-\overline\eta_j}} 
            \\[\fsmallskip]
            &\hspace*{5pt}=\left(-1\right)^{m_1} \overline C_1 \, C_2 \, \e^{\I\alpha_2} \, 
            \zeta^{m_1+m_2} \, \e^{-\zeta (t_2+\overline \gamma_1- \gamma_2)}
            \smashoperator{\prod_{\xi_j \in \Xi_1\setminus\Lambda_1}} \,
            \Bigl( 1 - \tfrac{\zeta}{-\overline\xi_j} \Bigr) \, 
            \e^{\frac{\zeta}{-\overline\xi_j}}
            \smashoperator{\prod_{\eta_j \in \Xi_2\setminus\Lambda_2}} \,
            \Bigl( 1 - \tfrac{\zeta}{\eta_j} \Bigr) \, \e^{\frac{\zeta}{\eta_j}} 
        \end{aligned}
    \end{equation}
    holds for every $\zeta$ in the complex plane.  Since the possibly infinite products
    above are again entire functions by Borel's theorem \cite{Mar77},
    the zeros on both sides of the equality have to coincide.  However, due to the
    assumption that the zeros $\xi_j$ and $\eta_j$ of the Laplace transforms $F$ and $H$
    are pairwise distinct, the zero sets $\Xi_1 \setminus \Lambda_1$ and
    $\Xi_2 \setminus \Lambda_2$ of $F_2$ and $H_2$ have to be invariant under reflection
    at the imaginary axis.

    Based on this observation, we can immediately conclude that the entire functions $F_2$
    and $H_2$ are invariant under reflection and conjugation up to an additional rotation
    and modulation.  More precisely, we obtain the identities
    \begin{align*}
      \overline{F_2\mleft( - \overline \zeta \mright)}
      &= \left( -1 \right)^{m_1} \overline C_1 \, \zeta^{m_1} \, \e^{-\zeta \overline
        \gamma_1}
        \smashoperator{\prod_{\xi_j \in \Xi_1\setminus\Lambda_1}} \,
        \Bigl( 1 - \tfrac{\xi}{-\overline \xi_j} \Bigr) \, \e^{\frac{\zeta}{- \overline \xi_j}}
      \\[\fskip]
      &= \left(-1\right)^{m_1} \e^{-2\I \arg C_1} \, \e^{-2 \zeta \Re [\gamma_1]} \, F_2
        \mleft( \zeta \mright)
        \intertext{and similarly}
        \overline{H_2\mleft( - \overline \zeta \mright)}
      &= \left(-1\right)^{m_2} \e^{-2\I \arg C_2} \, \e^{-2 \zeta \Re [\gamma_2]} \, H_2
        \mleft( \zeta \mright).
    \end{align*}
    Incorporating these identities in the representation of $\breve F$ and $\breve H$, we
    can describe the second solution pair in the frequency domain by
    \begin{align*}
      \breve F \mleft( \zeta \mright)
      &= \left(-1\right)^{m_1} \e^{\I(\alpha_1 + 2 \arg C_1)} \, \e^{-\zeta(t_1 - 2
        \Re[\gamma_1])} \, \overline{F \mleft( - \overline \zeta \mright)}
        \shortintertext{and}
        \breve H \mleft( \zeta \mright)
      &= \left(-1\right)^{m_2} \e^{\I(\alpha_2 + 2 \arg C_2)} \, \e^{-\zeta(t_2 - 2
        \Re[\gamma_2])} \, \overline{H \mleft( - \overline \zeta \mright)}.
    \end{align*}
    Hence, the continuous-time signals $\breve f$ and $\breve h$ are merely rotations and
    shifts of the reflected and conjugated signals $f$ and $h$.

    It remains to prove that the occurring rotations and shifts coincide.  For this
    purpose, we revisit equation \eqref{eq:cont-unknown-ref:sec-fact}.  Considering that
    the zeros and hence the possibly infinite products on both sides are equal, we can reduce
    \eqref{eq:cont-unknown-ref:sec-fact} to
    \begin{equation*}
        \left(-1\right)^{m_2} C_1 \, \overline C_2 \, \e^{\I \alpha_1} \, \e^{- \zeta (t_1
            - \gamma_1 + \overline \gamma_2)}
        = \left(-1\right)^{m_1} \overline C_1 \, C_2 \, \e^{\I \alpha_2} \, \e^{- \zeta (t_2
            + \overline \gamma_1 - \gamma_2)}
    \end{equation*}
    or, by rearranging and combining the individual factors, to
    \begin{equation*}
        \left(-1\right)^{m_1} \e^{\I \alpha_1+2 \arg C_1} \,
        \e^{-\zeta( t_1 - 2 \Re [\gamma_1])}
        = \left(-1\right)^{m_2} \e^{\I \alpha_2+2 \arg C_2} \,
        \e^{-\zeta( t_2 - 2 \Re [\gamma_2])},
    \end{equation*}
    which verifies our conjecture that the second solution pair $\breve f$ and $\breve h$
    coincides with the first solution pair $f$ and $h$ up to common trivial ambiguities.

    For the second case, where the first factor of \eqref{eq:cont-unknown-ref:fact-intf} is
    constantly zero, an analogous and slightly simpler argumentation yields the
    representations
    \begin{equation*}
        \breve F \mleft( \zeta \mright)
        = \e^{\I \alpha_1} \, \e^{- \zeta t_1} \, F \mleft( \zeta \mright)
        \quad\text{and}\quad
        \breve H \mleft( \zeta \mright)
        = \e^{\I \alpha_2} \, \e^{- \zeta t_2} \, H \mleft( \zeta \mright),
    \end{equation*}
    where the occurring rotations and time shifts again coincide.  
    \qed%
\end{Proof}

Our last approach to achieve the uniqueness of the continuous-time phase retrieval problem
considered in this section is again the idea of using interference measurements of the
unknown signal with a modulated version of the signal itself.  Generalizing the main
results of section~\ref{sec:interf-modul-sig}, we will establish two different theorems, which
show that each continuous-time signal in $L_2$ with compact support can be uniquely
recovered from an appropriate set of interference measurements.

\begin{Theorem}
    \label{the:cont-interf-polar-id}%
    Let $f$ be a continuous-time signal in $L^2$ with compact support.  Then the signal
    $f$ can be uniquely recovered up to rotation from its Fourier intensity
    $\absn{\Fourier [f]}$ and the interference measurements
    \begin{equation*}
        \abs{\Fourier\left[ f + \zeta_K^{-k} \, \e^{\I \mu \cdot} \, f \right]}
        \qquad \left( k =0, \dots, K-1; \mu \in M \right)
    \end{equation*}
    for an integer $K>2$ and every open neighbourhood $M$ around zero.%
\end{Theorem}

\begin{Proof}
    Due to the assumption that the unknown signal $f$ is a square-integrable function with
    compact support, the theorem of Paley-Wiener implies that the Fourier transform
    $\Fourier [f]$ is the restriction of an entire function and thus has to be continuous.
    Consequently, if the signal $f$ does not vanish everywhere, we can find a point
    $\omega_0$ together with an open neighbourhood where the Fourier transform
    $\Fourier [f]$ does not vanish.

    Similarly to the discrete-time version in Theorem~\ref{the:interf-polar-id}, the key element
    of the proof is to exploit the additional interference measurements
    $\absn{\fourier f (\cdot ) + \zeta_K^{-k} \, \fourier f ( \cdot - \mu)}$ by using the
    polarization identity in Lemma~\ref{lem:polar-id}.  We obtain
    \begin{equation*}
        \frac{1}{K} \sum_{k=0}^{K-1} \zeta_K^k \abs{\fourier f \mleft(\omega_0 \mright) +
            \zeta_K^{-k} \, \fourier f \mleft( \omega_0 - \mu \mright)}^2
        = \overline{\fourier f \mleft( \omega_0 \mright)} \, \fourier f \mleft( \omega_0 -
        \mu \mright)
    \end{equation*}
    for every $\mu$ in the open set $M$.  Writing the Fourier transform $\Fourier[f]$ in
    its polar representation $\absn{\Fourier[f]} \, \e^{\I\phi}$, we can now extract the
    relative phases $\phi(\omega_0-\mu) - \phi(\omega_0)$ from
    \begin{equation*}
        \frac{1}{K} \sum_{k=0}^{K-1} \zeta_K^k \abs{\fourier f \mleft(\omega_0 \mright) +
            \zeta_K^{-k} \, \fourier f \mleft( \omega_0 - \mu \mright)}^2
        = \abs{\fourier f \mleft( \omega_0 \mright)} \abs{ \fourier f \mleft( \omega_0 -
            \mu \mright)} \e^{\I(\phi(\omega_0-\mu) - \phi(\omega_0))}.
    \end{equation*}

    Like for the discrete-time counterpart, the considered phase retrieval problem can
    merely be solved up to rotations.  This enables us to define one phase
    $\phi(\omega_0)$ in the frequency domain arbitrarily.  Beginning from this initial
    phase, we can further determine the complete phase function $\phi$ and hence the
    Fourier transform $\Fourier [f]$ in a small open interval around $\omega_0$ by using
    the extracted relative phases.  Since the unknown Fourier transform $\Fourier [f]$ is
    the restriction of an entire function as discussed above, the unknown function
    $\Fourier [f]$ can be uniquely extended from the small interval to the complete
    frequency domain.  Using the inverse Fourier transform, we finally obtain the desired
    signal $f$.  \qed%
\end{Proof}

\begin{Theorem}
    \label{the:cont-interf-two-meas}%
    Let $f$ be a continuous-time signal in $L^2$ with compact support.  Then the signal
    $f$ can be uniquely recovered up to a rotation from its Fourier intensity
    $\absn{\Fourier [f]}$ and the interference measurements
    \begin{equation*}
        \abs{\Fourier\left[ f + \e^{\I\alpha_1} \, \e^{\I \mu \cdot} \, f \right]}
        \quad\text{and}\quad
        \abs{\Fourier\left[ f + \e^{\I\alpha_2} \, \e^{\I \mu \cdot} \, f \right]}
        \qquad \left(\mu \in M \right)
    \end{equation*}
    where $\alpha_1$ and $\alpha_2$ are two real numbers satisfying
    $\alpha_1 - \alpha_2 \ne \pi k$ for all integers $k$, and where $M$ is an open
    neighbourhood around zero.%
\end{Theorem}

\begin{Proof}
    Again, the crucial point to verify the assertion is the extraction of the relative
    phase from the given interference measurements.  Letting $\phi$ be the phase function
    of the unknown Fourier transform $\Fourier [f]$, and following the lines in the proof
    of the discrete-time counterpart (Theorem~\ref{the:two-interf-arb-rot}), we can determine
    the values
    \begin{equation*}
        \Re \mleft[ \e^{\I(\phi(\omega - \mu) - \phi(\omega)+ \alpha_1)} \mright]
        \quad\text{and}\quad
        \Re \mleft[ \e^{\I(\phi(\omega - \mu) - \phi(\omega)+ \alpha_2)} \mright]
        \addmathskip
    \end{equation*}
    and further the relative phase $\phi(\omega-\mu)-\phi(\omega)$ whenever
    $\Fourier[f](\omega)$ and $\Fourier[f](\omega-\mu)$ are non-zero by solving a linear
    equation system.  Based on the extracted relative phases
    $\phi(\omega-\mu)-\phi(\omega)$, we now can recover the unknown Fourier transform
    $\Fourier[f]$ as discussed in the previous proof of Theorem~\ref{the:cont-interf-polar-id}.
    \qed%
\end{Proof}

\begin{Remark}
    Since a holomorphic function is completely determined by the function values on an
    arbitrary set that has an accumulation point, we can relax the requirements on $M$ in
    Theorem~\ref{the:cont-interf-polar-id} and Theorem~\ref{the:cont-interf-two-meas}.  In fact, we
    can replace the condition that $M$ is an open neighbourhood around zero by the
    assumption that $M$ possesses at least one accumulation point.  \qed
\end{Remark}

\section*{Acknowledgements}

I gratefully acknowledge the funding of this work by the DFG in the frame work of the
SFB~755 \bq{Nanoscale photonic imaging} and of the GRK~2088 \bq{Discovering structure in
    complex data: Statistics meets Optimization and Inverse Problems.}

\bibliographystyle{alphadinUK}
{\footnotesize \bibliography{LITERATURE}}

\end{document}